\documentclass[12pt]{article}
\usepackage{amssymb}
\usepackage{graphicx}
\usepackage{amsfonts}
\usepackage{latexsym}
\usepackage{amsthm}
\usepackage{amsmath}

\usepackage{float}
\usepackage{subfigure}
\usepackage{bbm}
\usepackage{xcolor}

\usepackage{amssymb, amscd}

\usepackage{amsmath,amsfonts,graphicx,subfigure,url,amscd}

\usepackage{tikz}
\usetikzlibrary{decorations.pathreplacing,decorations.markings}
\usetikzlibrary{arrows}
\usepackage{mathtools}

\usepackage{url}
\usepackage[T1]{fontenc}

\usepackage{color}
\usepackage{tabu}

\numberwithin{equation}{section}

\newtheorem{thm}[equation]{Theorem}
\newtheorem{rem}[equation]{Remark}

\newtheorem{defin}[equation]{Definition}
\newtheorem{prop}[equation]{Proposition}

\newtheorem{lema}[equation]{Lemma}
\newtheorem{exam}[equation]{Example}

  \usetikzlibrary{shapes,arrows.meta,decorations.markings,calc}
  \tikzset{->-/.style={
    decoration = {
      markings,
      mark = at position .5 with {\arrow[scale=1.2]{latex}}
    },
    postaction = {decorate}
  }}
\tikzset{mid arrow/.style={postaction={decorate,decoration={
        markings,
        mark=at position .5 with {\arrow[#1]{stealth}}
      }}}
      }
\begin{document}
\date{May 1, 2020}
\title{{Generalized Tonnetz and discrete Abel-Jacobi map}}

\author{{Filip D. Jevti\'{c}} \\{\small Mathematical Institute}\\[-2mm] {\small SASA,   Belgrade}
\and Rade T. \v Zivaljevi\' c\\ {\small Mathematical Institute}\\[-2mm] {\small SASA,   Belgrade}}

\maketitle
\begin{abstract}\noindent
Motivated by classical Euler's {\em Tonnetz}, we introduce and study the combinatorics and topology of  more general simplicial complexes $Tonn^{n,k}(L)$ of {\em Tonnetz type}. Out main result is that for a sufficiently generic choice of parameters the generalized tonnetz $Tonn^{n,k}(L)$ is a triangulation of a $(k-1)$-dimensional torus $T^{k-1}$. In the proof we construct and use the properties of a {\em discrete Abel-Jacobi map}, which takes values in the torus $T^{k-1} \cong \mathbb{R}^{k-1}/\Lambda$ where $\Lambda \cong \mathbb{A}^\ast_{k-1}$ is the permutohedral lattice.
\end{abstract}



\renewcommand{\thefootnote}{}
\footnotetext{This work was supported by the Serbian Ministry of Education, Science and Technological Development through Mathematical Institute of the Serbian
Academy of Sciences and Arts.}

\section{Introduction}
\label{sec:intro}

In his seminal work on music theory \textit{``Tentamen novae theoriae musicae ex certissismis harmoniae principiis dilucide expositae''} (1739), Leonhard Euler introduced a lattice diagram -- \textit{Tonnetz} -- representing the classical tonal space.
In more recent interpretations this diagram is identified as a triangulation of a torus with 24 triangles representing all the major and minor chords.
If the equal tempered scale is identified with $\mathbb{Z}_{12}$, the Tonnetz can be described as
\begin{align*}
    Tonnetz=\left\lbrace \{x,x+3,x+7\} \mid x \in \mathbb{Z}_{12} \right\rbrace \cup \left\lbrace \{x,x+4,x+7\} \mid x \in \mathbb{Z}_{12} \right\rbrace.
\end{align*}
Notice that if $\{x,y,z\} \in Tonnetz$ then $\{x-y,y-z,z-x\} = \pm\{ 3, 4, 5\}$.
This serves as an inspiration to introduce and study  more general complexes of ``Tonnetz type.''

\subsection{Generalized Tonnetz}
\label{sec:intro-intro}

Suppose that $L = \{l_i\}_{i=1}^k$ is a collection of $k$ positive integers which add up to $n$
\[
   l_1 + l_2 + \dots + l_k = n \, .
\]
We say that a collection $L$ is generic if for each pair $I, J\in 2^{[n]}$ of subsets of $[n]$
\begin{equation}\label{eqn:generic}
     \sum_{i\in I} l_i = \sum_{j\in J} l_j  \, \Rightarrow   \,  I = J \, .
\end{equation}
A collection $L$ is reduced if the largest common divisor of all $l_i$ is $1$,
\begin{equation}\label{eqn:divisor}
 \langle l_1, l_2,\dots, l_k \rangle   = 1 \, .
\end{equation}

\medskip\noindent
{\bf Caveat:} From here on we  identify elements of $[n]$ with the corresponding elements (congruence classes) in  the additive group $\mathbb{Z}_n = \{0,1,\dots, n-1\}$ of
integers, modulo $n$. A standard geometric model for this set is $\mathbb{V}_n = \{\epsilon  \mid \epsilon^n =1  \}$,  the set of vertices of a regular $n$-gon. For this reason we in principle assume a counterclockwise orientation on the unit circle $S^1$ where the $n$-gon is inscribed. However for some readers it may be more natural to use (occasionally) more traditional presentation of the (classical) Tonnetz, with clockwise orientation and $n=12$ occupying its ``usual'' place.

\medskip
Each ordered pair $(x,y)$ of elements in $\mathbb{Z}_n$ (respectively $[n]$ or  $\mathbb{V}_n$) defines an interval $I_{x,y} = \{x,x+1,\dots, y\}\subset \mathbb{Z}_n$. The length of this interval  is $\mathfrak{L}(I_{x,y}) = \vert I_{x,y}\vert - 1 = y-x \in \mathbb{Z}_n$.

\medskip
If a set $\tau = \{v_0, v_1, \dots, v_t\}$ is a subset of $\mathbb{Z}_n$ we always assume the cyclic order $v_0 \prec v_1 \prec\dots\prec v_t$ of its elements, i.e.
$v_j - v_i\in \{1,2,\dots, n-1\}$ for each pair $i<j$ of indices.

\begin{defin}\label{def:Tonn}
A generalized tonnetz $Tonn^{n,k}(L)\subseteq 2^{[n]}$
is a $(k-1)$-dimensional simplicial complex whose maximal simplices are
\begin{align}\label{align:tonnetz-max-simplex}
    \Delta(x;\sigma)=\{x,x+l_{\sigma(1)},x+l_{\sigma(1)}+l_{\sigma(2)},\ldots,x+l_{\sigma(1)}+\ldots+l_{\sigma(k-1)}\}
\end{align}
where $x\in \mathbb{Z}_n$ and $\sigma\in \Sigma_k$ is a permutation.
\end{defin}

Our main result is the following theorem which claims that a (sufficiently generic) generalized Tonnetz is also a triangulation of a torus, as its classical counterpart.

\begin{thm}\label{thm:main}
Suppose that $L = \{l_i\}_{i=1}^k$ is generic and reduced in the sense of (\ref{eqn:generic}) and (\ref{eqn:divisor}). Then the generalized tonnetz $Tonn^{n,k}(L)$ is a triangulation of a $(k-1)$-dimensional
torus $T^{k-1} := (S^1)^{k-1}$.
\end{thm}

The central idea in the proof of Theorem \ref{thm:main} is to identify the (triangulated)  torus $T^{k-1}\cong\mathbb{R}^{k-1}/\Lambda$ as a {\em combinatorial Jacobian}, i.e. as the target of a (discrete) Abel-Jacobi map   $J :   Tonn^{n,k}(L) \longrightarrow \mathbb{R}^{k-1}/\Lambda$ where $\Lambda \cong \mathbb{A}^\ast_{k-1}$ is a permutohedral lattice.

\subsection{Discrete Abel-Jacobi map}
\label{sec:abel-jacobi}

The classical Abel-Jacobi map is a map from an algebraic curve $S$ (Riemann surface of genus $g$) into the torus $\mathbb{C}^g/\Lambda$ where $\Lambda\subset \mathbb{C}^g$ is the lattice of periods. More explicitly there exist $g$ linearly independent holomorphic differentials $\omega_1,\dots,\omega_g$ on $S$ and if $\{c_j\}_{j=1}^{2g}\subset H_1(S; \mathbb{Z})$ is a collection of basic cycles then the vectors
\[
        v_j = \langle c_j, \omega \rangle = (\int_{c_j} \omega_1, \dots, \int_{c_j} \omega_g)
\]
form a basis of a lattice $\Lambda$.  Then the Abel-Jacobi map $J : S \rightarrow \mathbb{C}^g/\Lambda $ is defined by
\begin{equation}\label{eqn:Abel-Jac-holomorphic}
    J(p) = (\int_{p_0}^p \omega_1, \dots, \int_{p_0}^p \omega_g) \, \mbox{ {\rm mod} } \Lambda \, .
\end{equation}

In the special case when the curve $S$ is elliptic the Jacobi map is an isomorphism.

\medskip
In  analogy with this construction, we describe explicit simplicial cocycles $\omega_i \, (1\leqslant i \leqslant k)$ on $Tonn^{n,k}(L)$,  which play the role of holomorphic differentials and allow us to construct the corresponding ``discrete Abel-Jacobi map.'' This can be  compared to the use of discrete Abel-Jacobi maps in the construction of {\em standard realizations of maximal abelian covers of graphs} in topological crystallography, see \cite{Su08, Su12, Su13} and \cite{Ba}.

\section{$Tonn^{n,k}(L)$ is a manifold}
\label{sec:Is-manifold}

We begin our analysis of complexes of  Tonnetz type by showing that the irreducibility condition (\ref{eqn:divisor}) can always be assumed,  without an essential loss of generality.

\begin{prop}\label{prop:reduction}
The (geometric realization of the) complex $Tonn^{pn,k}(pL) \subseteq 2^{[pn]}$, where $pL := \{pl_1, pl_2,\dots, pl_k\}$, is homeomorphic to the disjoint union of $p$ copies of $Tonn^{n,k}(L)$.
\end{prop}

\medskip\noindent
{\bf Proof:} As a consequence of (\ref{align:tonnetz-max-simplex}) if $\sigma = \{v_1, v_2,\dots, v_k\}\in Tonn^{pn,k}(pL)$ then $v_i \equiv v _j\, (\mbox{\rm mod } p)$ for each $i,j\in [k]$.
It follows that $Tonn^{pn,k}(pL)$ is a disjoint union of its $p$ subcomplexes  $T_j \cong Tonn^{n,k}(L)$ where $T_j$ is spanned by  vertices in the same $\mathbb{Z}_n$-coset of the group $\mathbb{Z}_{pn}$.
\hfill $\square$

\medskip
Unlike the irreducibility condition (\ref{eqn:divisor}), the genericity condition (\ref{eqn:generic}) is essential for  the proof of the following proposition.

\begin{prop}\label{prop:manifold}
  $Tonn^{n,k}(L)$ is a connected, combinatorial manifold if  the ``length vector'' $L = (l_1,\dots, l_k)$ is both generic and reduced, in the sense of (\ref{eqn:generic}) and (\ref{eqn:divisor}).
  Moreover, the links of vertices are isomorphic to  boundaries of simplicial polytopes dual to $(k-1)$-dimensional permutohedra.
\end{prop}

\medskip\noindent
{\bf Proof:}  By definition $\tau = \{v_0,v_1,\dots, v_t\}\in Tonn^{n,k}(L)$ if and only there is a partition $[k] = I_0\sqcup I_1\sqcup\dots\sqcup I_t$  such that for each $i$  the length of the interval $I_{v_{i}, v_{i+1}} \, (v_{t+1}:= v_0)$ is
\[
    \mathfrak{L}(I_{v_{i}, v_{i+1}}) = \sum_{j\in I_i} l_j\, .
\]

For a chosen vertex $v = v_0$ the face poset of the star $Star(v_0)\subseteq Tonn^{n,k}(L)$ is isomorphic to the poset of all ordered partitions of $[k]$ (here we use the genericity of $L$), where  the top dimensional simplices in $Star(v_0)$ correspond to the finest ordered partitions of $[k]$.
Note that the finest ordered partitions of $[k]$ are in $1$-$1$ correspondence with permutations of $[k]$.

Recall \cite[Example 0.10]{Ziegler-book} that the face poset of a $(k-1)$-dimensional permutohedron $Perm_{k-1}$ is
also the poset of all ordered partitions of $[k]$, but with the reversed ordering.
(The finest partitions/permutations correspond to the vertices of the permutohedron.)
 It immediately follows that $Link(v_0) \cong \partial Q_{k-1}$ where $Q_{k-1} = Perm_{k-1}^\circ$ is the (simplicial) polytope polar to the permutohedron.

\medskip
More generally, the link $Link(\tau)$ of $\tau = \{v_0,v_1,\dots, v_t\}$ is isomorphic
to the join
\[
     Link(\tau) = \partial Q_{s_0-1} \ast \dots \ast \partial Q_{s_t-1}
\]
where  $s_j = \vert I_j\vert$ is the cardinality of the set $I_j$.
Consequently, $Tonn^{n,k}(L)$ is indeed a manifold. To show that it is connected, it is sufficient to show that consecutive vertices $x$ and $x+1$ are connected. Indeed, since $L$ is reduced we obtain the relation
\begin{align*}
a_1l_1+\ldots + a_kl_k=1
\end{align*}
for some $a_1,\ldots,a_k \in \mathbb{Z}$, which describes a sequence of edges connecting $x$ and $x+1$. \hfill $\square$

\begin{rem}\label{rem:zonotope}{\rm
The following geometric model for the complex $Star(v_0)$ can be used for an alternative proof of Proposition \ref{prop:manifold}. Let $c_i\in\mathbb{R}^{k-1}\, (i=1,\dots, k)$ be a spanning set of vectors such that $c_1+c_2+\dots+ c_k =0$.  Let $Z = [0, c_1]+\dots+ [0,c_k]\subset \mathbb{R}^{k-1}$ be the Minkowski sum of line segments $I_j = [0, c_j]$. Then the zonotope $Z$ admits a triangulation where the maximal simplices $\Sigma_\pi \, (\pi\in S_k)$,  indexed by permutations, are the following
\begin{equation}\label{eqn:zono}
    \Sigma_\pi = {\rm Conv}\{c_{\pi(1)}, c_{\pi(1)} + c_{\pi(2)}, \dots,  c_{\pi(1)}+\dots + c_{\pi(k)}\} \, .
\end{equation} }
\end{rem}
This triangulation of $Z$ is isomorphic to $Star(v_0)$ which can be proved by comparing (\ref{align:tonnetz-max-simplex}) and (\ref{eqn:zono}).

\medskip
A very special case of Theorem \ref{thm:main} can be established by an elementary, direct argument.

\begin{prop}\label{prop:k=3}
 Let $L = (l_1, l_2, l_3)$ be a reduced, generic length vector. Then the associated, $2$-dimensional generalized Tonnetz  $Tonn^{n,3}(L)$ is a triangulation of the $2$-dimensional torus $T^2 = (S^1)^2$.
\end{prop}

\medskip\noindent
{\bf Proof:}  In light of Propositions \ref{prop:reduction} and \ref{prop:manifold} $Tonn^{n,3}(L)$ is a connected $2$-manifold.
Assume $l_1< l_2 < l_3$.  It is not difficult to see that the $f$-vector of  $T = Tonn^{n,3}(L)$ is $f(T)=(n,3n,2n)$, hence $\chi(T)=0$.

\medskip
The complex $Tonn^{n,3}(L)$ is orientable. Indeed, all triangles in $Tonn^{n,3}(L)$ fall into two classes. Generalized ``major triads'' are the triangles $\tau = \{v_0, v_1, v_2 \}$  where $v_1 - v_0 = l_1, v_2-v_1 = l_2$ and $v_0- v_2 = l_3$ (for some circular order of vertices of $\tau$). These simplices are positively oriented. Negatively oriented are generalized ``minor triads'', i.e.  the triangles $\tau = \{v_0, v_1, v_2 \}$  where $v_1 - v_0 = l_1, v_2-v_1 = l_3$ and $v_0- v_2 = l_2$.

\medskip
Summarizing, $Tonn^{n,3}(L)$ is a connected, orientable $2$-manifold with vanishing Euler characteristic, hence it must be the torus $T^2$.  \hfill $\square$

\subsection{Euler characteristic and $f$-vector of $Tonn^{n,k}(L)$}

\begin{prop}
Suppose that $L$ is generic and reduced and let  $f(T) = (f_0, f_1,\dots, f_{k-1})$ be the $f$-vector of the generalized Tonnetz $T = Tonn^{n,k}(L)$. Then \begin{align}\label{eqn:Stirling}
    f_{m-1}=n\frac{P(k,m)}{m} =n\frac{m!}{m}S(k,m)
\end{align}
where $P(k,m)$ is number of ordered partitions partition $I_{1} \sqcup \ldots \sqcup I_{m} = [k]$ and $S(k,m)$ is a Stirling number of the second kind.
\end{prop}

  \medskip\noindent
{\bf Proof:} If  $S$ is a $(m-1)$-dimensional face of $Tonn^{n,k}(L)$ then  there exists $x\in \mathbb{Z}_{n}$ and a partition $[n]=I_1\sqcup \ldots \sqcup I_m$ such that
\begin{align*}
    S=\left\lbrace x,x+\mu_{L}(I_1),x+\mu_{L}(I_1\sqcup I_2),\ldots,x+\mu_{L}(I_1\sqcup \ldots \sqcup I_{m-1})\right\rbrace,
\end{align*}
where $\mu_{L}(I)=\sum_{i \in I} l_i$.
The first equality in the formula (\ref{eqn:Stirling}) is an immediate consequence.
(Since $L$ is generic if $\mu_L(I) = \mu_L(J)$ then $I=J$.)
The second follows from the equality $P(k,m)=m!S(k,m)$, where $S(k,m)$ is a Stirling number of the second kind.   \hfill $\square$

\begin{prop}
  If the vector $L$ is generic and reduced then the Euler characteristic of $Tonn^{n,k}(L)$ is 0.
\end{prop}

\medskip\noindent
{\bf Proof:} Let $T=Tonn^{n,k}(L)$. Its Euler characteristic is
\begin{align*}
    \chi(T) =n \sum_{m=1}^{k} \frac{(-1)^{m+1}m!}{m}S(k,m).
\end{align*}

Using the well-known recurrence for Stirling numbers
\begin{align*}
    S(k,m)=mS(k-1,m)+S(k-1,m-1)
\end{align*}
we obtain
\begin{align*}
    \chi(T) &= n \left( \sum_{m=1}^k \frac{(-1)^{m+1}m!}{m}mS(k-1,m) + \sum_{m=1}^k \frac{(-1)^{m+1}m!}{m}S(k-1,m-1) \right)\\
    &=n \left( \sum_{m=1}^k (-1)^{m+1}m!S(k-1,m) + \sum_{q=0}^{k-1} (-1)^{q}q!S(k-1,q) \right)\\
    &=n \left( (-1)^{k+1}k!S(k-1,k) +S(k-1,0) \right)\\
    &= 0
\end{align*}
\hfill $\square$

\subsection{Fundamental group of $Tonn^{n,k}(L)$}

A consequence of Theorem \ref{thm:main} is that the fundamental group of a generalized Tonnetz $Tonn^{n,k}(L)$ is free abelian of rank $k-1$, provided the vector  $L$ is generic and reduced. Proposition \ref{prop:abelian}  is a key step in the direction of this result. Before we commence the proof, let us make some general observation about the edge-path groupoid of the Tonnetz $Tonn^{n,k}(L)$.

\medskip
Each edge-path connecting vertices $a = v_0$ and $b=v_m$ is of the form $\alpha = X_1X_2\cdots X_m$ where $X_i = \overrightarrow{v_{i-1}v_i}$ is an oriented edge ($1$-simplex) in $Tonn^{n,k}(L)$.

\medskip
Recall (Section \ref{sec:intro-intro}) that $I_{X} = I_{u,v}\subset \mathbb{Z}_n$ is the (oriented) interval, corresponding to $X = \overrightarrow{u v}$. (With a slight abuse of language we use the same notation for the corresponding arc in $S^1$.)

Let the $L$-type $X^L$ of $X$  be defined as the unique non-empty subset $I\subset [k]$ such that $\mathfrak{L}(I_{u, v}) = \sum_{j\in I} l_j$.

\medskip
We say that $X = \overrightarrow{u v}$ is {\em atomic} if either $v = u + l_i$ or $u = v+l_i$ for some $i\in [k]$.  If $X$ is positively oriented, i.e. if $v = u + l_i$, then we call it $\oplus$-atomic (similarly  $\circleddash$-atomic if $u = v + l_i$).

\medskip
Note that the $L$-type of an $\oplus$-atomic oriented $1$-simplex $X = \overrightarrow{u v}$ is a singleton $X^L = \{i\}$ (we say that $X$ is of type $i$), while the $L$-type of the associated $\circleddash$-atomic $1$-simplex $X^{-1} = \overrightarrow{v u}$ is $[k]\setminus \{i\}$.

\medskip
The following lemma is an immediate consequence of Definition \ref{def:Tonn}.

\begin{lema}\label{lema:atom}
  Each oriented $1$-simplex $X = \overrightarrow{u v}$ is homotopic $X \cong Y_1Y_2\cdots Y_t$ (relative to the end-points $u$ and $v$) to a product of $\oplus$-atomic $1$-simplices $Y_j$. Moreover,  one can read off the $L$-type of $X$ from this representation as,
  \begin{equation}\label{eqn:atom-decomposition}
  X^L = \{ Y_1^L, Y_2^L, \dots , Y_t^L\} \, .
 \end{equation}
\end{lema}

\begin{figure}[htb]
	\centering
	\begin{tikzpicture}[scale=1]
	\newcommand{\ang}{360/12}	% angle
	\newcommand{\R}{3}	% radius
	\draw (0,0) circle(3);
        
    \foreach \x in {0,...,8,9,10,11} {
     	\draw ({\x*\ang}:{\R-.1}) -- ({\x*\ang}:{\R+.1});
     	%\node () at ({\x*\ang}:{\R+.5}) {\x};
  	}
        
    \draw[very thick,->-] ({4*\ang}:\R) to node[above right] {$X$} ({11*\ang}:{\R}) ;
    \draw[very thick,->-] ({4*\ang}:\R) to node[below right] {$Y_1$} ({6*\ang}:{\R}) ;
    \draw[very thick,->-] ({6*\ang}:\R) to node[above right] {$\ldots$} ({8*\ang}:{\R}) ;
    \draw[very thick,->-] ({8*\ang}:\R) to node[above] {$Y_t$} ({11*\ang}:{\R}) ;
\end{tikzpicture}
    \caption{$X=Y_1\ldots Y_t$}
    \label{fig-1}
\end{figure}
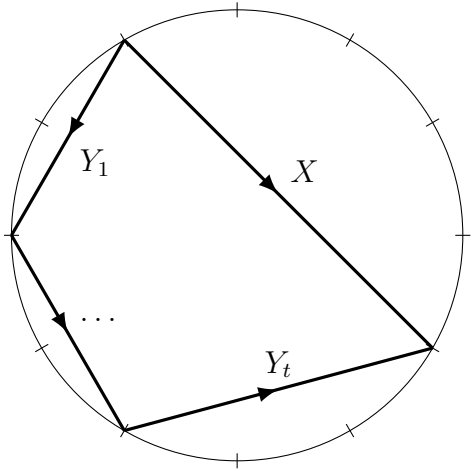

The following lemma (see Figure \ref{fig-2} for a visual proof) shows that  we can rearrange and group $\oplus$-atomic $1$-simplices according to their type.

\begin{lema}\label{lema:+atom}
An edge-path which is a product $XY$ of two $\oplus$-atomic $1$-simplices $X$ and $Y$, respectively of type $i$ and $j$ (where $i\neq j$) is homotopic (rel the end points) to a product $Y'X'$ of two $\oplus$-atomic $1$-simplices, where $type(Y') = j$ and $type(X') = i$.
\end{lema}

\begin{figure}[htb]
	\centering
	\begin{tikzpicture}[scale=1]
	\newcommand{\ang}{360/12}	% angle
	\newcommand{\R}{3}	% radius
	\draw (0,0) circle(3);
        
    \foreach \x in {0,...,8,9,10,11} {
     	\draw ({\x*\ang}:{\R-.1}) -- ({\x*\ang}:{\R+.1});
     	%\node () at ({\x*\ang}:{\R+.5}) {\x};
  	}

    \draw[very thick] ({3*\ang}:\R) to  ({8*\ang}:{\R}) ;
    \draw[very thick,->-] ({8*\ang}:\R) to node[above] {$X$} ({11*\ang}:{\R}) ;
    \draw[very thick,->-] ({11*\ang}:\R) to node[right] {$Y$} ({3*\ang}:{\R}) ;
    \draw[very thick,->-] ({8*\ang}:\R) to node[above] {$Y$} ({0*\ang}:{\R}) ;
    \draw[very thick,->-] ({0*\ang}:\R) to node[right] {$X$} ({3*\ang}:{\R}) ;
        
    \draw[very thick,->-] (5,-1.5) to node[left] {$X$} (5,1.5);
    \draw[very thick,->-] (9,-1.5) to node[right] {$X$} (9,1.5);
    \draw[very thick,->-] (5,1.5) to node[above] {$Y$} (9,1.5);
    \draw[very thick,->-] (5,-1.5) to node[below] {$Y$} (9,-1.5);
    \draw[very thick] (5,-1.5) to (9,1.5);
\end{tikzpicture}
    \caption{$XY= YX$}
    \label{fig-2}
\end{figure}

\begin{prop}\label{prop:abelian}
If the length vector $L$ is generic and reduced then the fundamental group $\pi_1(Tonn^{n,k}(L))$ of the generalized Tonnetz  is abelian.
\end{prop}

\medskip\noindent
{\bf Proof:} Suppose that $v_0$ is the chosen base point and assume that $\alpha$ and $\beta$ are two edge-paths (loops) based at $v_0$. We are supposed to show that the edge paths $\alpha\beta$ and $\beta\alpha$ are homotopic (rel $v_0$).

\medskip
By Lemma \ref{lema:atom} we are allowed to assume that both $\alpha = Y_1Y_2\cdots Y_s$ and $\beta = Z_1Z_2\cdots Z_t$ are  products of $\oplus$-atomic $1$-simplices.

Use Lemma \ref{lema:+atom} to rearrange atoms in the product $\alpha\beta$ (similarly $\beta\alpha$) and write it as a product $\alpha\beta = A_1A_2\dots A_k$, where $A_i$ is the product of $\oplus$-atoms of type $i$. (Some $A_i$ may be empty words.)

\medskip
Observe that (as a consequence of Lemma \ref{lema:+atom}) the length of the word $A_i$ is equal to the number of type $i$ $\oplus$-simplices in the product $\alpha\beta$.

\medskip
If $\beta\alpha = A'_1A'_2\dots A'_k$ is the corresponding regrouped presentation of $\beta\alpha$ we observe that $A_j = A'_j$ for each $j$. This completes the proof of the proposition.  \hfill $\square$

\section{Canonical cycles and cocycles in $Tonn^{n,k}(L)$}
\label{sec:canonical}

We already know (Proposition \ref{prop:manifold}) that a generalized Tonnetz $T = Tonn^{n,k}(L)$ is a connected complex.
Since, according to Proposition \ref{prop:abelian}, the fundamental group $\pi_1(T)$ is abelian, it is isomorphic to the first homology group $H_1(T; \mathbb{Z}) := Z_1/B_1$, where $Z_1$ and $B_1$ are the corresponding groups of cycles and boundaries.

\medskip
When working with the homology group it is more customary to use additive notation. For example the $\oplus$-atom decomposition   $X \cong Y_1Y_2\cdots Y_t$ from Lemma \ref{lema:+atom} can be rewritten as the following equality (in homology),
$X = \sum_{i=1}^t Y_i$.

\medskip
In this section the emphasis is on (co)homology so here we follow the additive notation.

\begin{defin}\label{def:basic-cocycles}  The cochains  $\theta_{i,j}\in C^1 = Hom(C_1; \mathbb{Z})$, where $1\leqslant i\neq j\leqslant k$, are  defined on $\oplus$-atomic $1$-simplices as follows:
\[
\theta_{i,j}(Y) = \left\{ \begin{array}{ll}
         +1 & \mbox{if $Y$ is of $L$-type $i$}\\
        -1 & \mbox{if $Y$ is of $L$-type $j$}\\
        0 & \mbox{if  the $L$-type of $Y$ is neither $i$ nor $j$. } \end{array} \right.
        \]
If $X = \overrightarrow{uv}$ is an oriented $1$-simplex and  $X \cong Y_1Y_2\cdots Y_t$ its $\oplus$-atom decomposition from Lemma \ref{lema:atom}, then by definition \begin{equation}\label{eqn:extension}
   \theta_{i,j}(X) = \sum_{m=1}^t \theta_{i,j}(Y_m)  \, .
\end{equation}
On other oriented $1$-chains they are extended by linearity.
\end{defin}

\begin{prop}\label{prop:theta}
The cochain $\theta_{i,j}$ is well defined. Moreover, it is a cocycle which defines an element of $H^1(T; \mathbb{Z})$. These classes (cocycles) are referred to as ``elementary classes'' defined on $T = Tonn^{n,k}(L)$.
\end{prop}

\medskip\noindent
{\bf Proof:} We check that $\theta_{i,j}$ is well defined by showing that possibly different ways to extend $\theta_{i,j}$ lead to the same result. Essentially the only case when this happens is when we evaluate $\theta_{i,j}(-X) = \theta_{i,j}(X^{-1})$, where
$X = \overrightarrow{uv}$ and $X^{-1} = \overrightarrow{vu}$ (by formula (\ref{eqn:extension})) expecting to obtain the result $-\theta_{i,j}(X)$.

This is indeed the case since $\sum_{m=1}^k \theta_{i,j}(Y_m) = 0$, where $Y_m$ is a $\oplus$-atomic $1$-simplex of type $m$, for each $m\in [k]$. Similarly we obtain that the coboundary
\[
     \delta\theta_{i,j}(\tau) = \theta_{i,j}(\overrightarrow{u_0 u_1}) + \theta_{i,j}(\overrightarrow{u_1 u_2}) + \theta_{i,j}(\overrightarrow{u_2 u_0}) = 0
\]
is zero for each (oriented) $2$-simplex $\tau = \{u_0, u_1, u_2\}$.  \hfill $\square$

\begin{defin}\label{def:omega}
For each $i\in [n]$ let $\omega_i$ be the cocycle defined by
$ \omega_i  = \sum_{j\neq i} \theta_{i,j} $.
More explicitly, $\omega_i$ is the unique $1$-cocycle defined on $Tonn^{n,k}(L)$ such that for each $\oplus$-atomic $1$-simplex $Y$
\[
\omega_i(Y) =
\left\{ \begin{array}{ll}
         k-1 & \mbox{if $Y$ is of type $i$}\\
        -1 & \mbox{if $Y$ is of type $j\neq i$.} \end{array} \right.
\]
These cocycles are referred to as the ``canonical'' cocycles defined on $Tonn^{n,k}(L)$.
\end{defin}

As an immediate consequence of the definition we obtain the following relation
\begin{equation}\label{eqn:omega-zavisnost}
\omega_1+\omega_2+\dots+ \omega_k  = 0 \, .
\end{equation}

The complex $Tonn^{n,k}(L)$ also has naturally defined  $1$-cycles.

\begin{defin}\label{def:cycles} For  $i\in [k]$ let $c_i$   be the $1$-cycle  defined by
 $c_i := \sum_{x\in [n]} E_x^i$
  where  $E_x^i = \overrightarrow{xy}$ (for $x\in [n]$ and $i\in [k]$) is the $\oplus$-atomic  $1$-simplex of type $i$ with end-points $x$ and $y=x +{l_i}$.
\end{defin}

\begin{prop}\label{prop:c-zavisnost} If  $[c_i]\in H_1(Tonn^{n,k}(L); \mathbb{Z})$ is the homology class of the cycle $c_i$ then
\begin{equation}\label{eqn:c-zavisnost}
[c_1]+[c_2]+\dots+ [c_k]  = 0 \, .
\end{equation}
\end{prop}

\medskip\noindent
{\bf Proof:} Informally, the cycle $c_i$ is the sum of all $\oplus$-atomic $1$-simplices of type $i$. They can be concatenated to form $d$ irreducible cycles of length $q$, where $n=qd$ and $d = (n,l_i)$. For each $x\in \mathbb{Z}_n$ the cycle $\mathbb{E}_x = E_x^1 + E_{x+l_1}^2+\dots+ E^k_{l_1+\dots+ l_{k-1}}$ is trivial by Definition \ref{def:Tonn}. Since $\sum_{i=1}^k c_k = \sum_{x\in \mathbb{Z}_n} \mathbb{E}_x$ the equality (\ref{eqn:c-zavisnost}) is an immediate consequence. \hfill $\square$

\medskip
The following proposition implies that aside from (\ref{eqn:omega-zavisnost}) and (\ref{eqn:c-zavisnost}) there are essentially no other relations among $\{\omega_i\}_{i\in [k]}$ and $\left\lbrace[c_i]\right\rbrace_{i\in [k]}$.

\begin{prop}\label{prop:dual-classes-matrix}
 Let $\langle \cdot, \cdot\rangle$ be the pairing between the cohomology and homology classes and let $M = [m_{i,j}]_{i,j=1}^{k-1}$ be a $(k-1)\times (k-1)$-matrix where $m_{i,j} = \langle\omega_i, c_j\rangle$. Then  ${\rm det}(M) = n(nk)^{k-2}$.
\end{prop}

\medskip\noindent
{\bf Proof:}  By direct calculation we have

\begin{equation}\label{eqn:gram-matrix}
  \langle \omega_i, c_j \rangle  = \sum_{\nu\neq i} \langle \theta_{i,\nu}, c_j \rangle
  =  \left\{ \begin{array}{ll}
         n(k-1) & \mbox{if $i = j$}\\
        -n & \mbox{if $i\neq j$} \end{array} \right.
\end{equation}

It follows that $M$ is a circulant matrix with the associated polynomial equal to $f(x) = n(k-1)-n(x+x^2+\dots+ x^{k-2})$. Recall that  the determinant of the circulant matrix with the associated polynomial $f(x)$ is equal to $\prod_{j=1}^{k-1} f(\epsilon_j)$, where $\epsilon_j$ are solutions of the equation $x^{k-1}-1=0$. From here it immediately follows that $det(M) = n(nk)^{k-2}$.

\hfill $\square$

\section{Homology $H_1(Tonn^{n,k}(L))$}
\label{sec:homology}

In this section we complete the analysis and summarize  our knowledge about the homology group $H_1(T; \mathbb{Z})$ of the generalized Tonnetz $T=Tonn^{n,k}(L)$.

\medskip
We already know (Section \ref{sec:Is-manifold}) that each homological $1$-cycle  has a (multiplicative) representation $X = Y_1Y_2\cdots Y_t$ where $Y_i$ are $\oplus$-atomic $1$-simplices. We also write $X  = {}_{v_0}Y_1Y_2\cdots Y_t$  when we want to emphasize that the initial vertex of $Y_1$ (playing the role of the base point of the loop $X$) is $v_0$.

  If $Y_j$ is of $L$-type $i_j$ then we can also (symbolically) record this information as the word
\begin{equation}\label{eqn:record-types}
    X  = E_{i_1}E_{i_2}\dots E_{i_t}  = {}_{v_0}E_{i_1}E_{i_2}\dots E_{i_t} \, ,
\end{equation}
where $E_j$ denotes a step of length $l_j$ in the positive direction.
Most of the time we can safely remove the base point $v_0$ from the notation.
For example the representation $X = E_2E_1E_1E_1E_3E_3 = E_2E_1^3E_3^2$ describes an edge-path which begins at $v_0$, makes one step of type $2$, then three steps of type $1$ and finally two steps of type $3$.

\medskip
Lemma \ref{lema:+atom} can be interpreted as a symbolic (edge-path) relation $E_iE_j= E_jE_i$ which says that one can interchange two consecutive steps (as in Fig. 2) without changing the homotopy type of the edge-path (rel end-points). From here we easily deduce the following proposition.

\begin{prop}\label{prop:monomial-cycles}
Each homological cycle X  has a representation $X = E_1^{p_1}E_2^{p_2}\cdots E_k^{p_k}$ where $p_1, \dots, p_k$ are non-negative integers such that $p_1l_1+p_2l_2+\dots + p_kl_k = p_0n$ for some $p_0\geqslant 0$.
Moreover, if $Y = E_1^{p'_1}E_2^{p'_2}\cdots E_k^{p'_k}$  has a similar representation, where  $p'_1l_1+p'_2l_2+\dots + p'_kl_k = p'_0n$, then the cycles $X$ and $Y$ are homologous if and only if
\begin{equation}\label{eqn:vector-difference}
  (p_1,\dots, p_k) - (p'_1,\dots, p'_k) \in \mathbb{Z}\, \mathbbm{1} \, .
\end{equation}
where $\mathbbm{1} = (1,\dots, 1)\in \mathbb{Z}^k$ and $\mathbb{Z} \, \mathbbm{1} = \{m\mathbbm{1} \vert\, m\in \mathbb{Z}\}$.
\end{prop}

\noindent
{\bf Proof:} If the relation (\ref{eqn:vector-difference}) is satisfied then $X$ and $Y$ are clearly homologous since $E = E_1E_2\dots E_k$ is a trivial cycle.

Conversely, suppose that $X$ and $Y$ are homologous. Then,
\[
   p_i - p_j = \theta_{i,j}(X) = \theta_{i,j}(Y) = p'_i - p'_j
\]
for each $i<j$ and the relation (\ref{eqn:vector-difference}) follows.  \hfill $\square$

\medskip
As an immediate consequence we obtain the following representation of the first homology group of the generalized Tonnetz $Tonn^{n,k}(L)$ as a lattice of rank $(k-1)$ in a hyperplane $H_0^L \subset \mathbb{R}^k$.

\begin{thm}\label{thm:main-homology}
Let $H_0^L = \{ x\in \mathbb{R}^k \mid \, \langle x, L \rangle = x_1l_1+\dots + x_kl_k = 0 \}$ be the central hyperplane in $\mathbb{R}^k$, orthogonal to $L$. Then there is an isomorphism
\begin{equation}
    H_1(Tonn^{n,k}(L); \mathbb{Z}) \longrightarrow  H_0^L\cap \mathbb{Z}^k
\end{equation}
where  $H_0^L\cap \mathbb{Z}^k$ is a free abelian group (lattice) of rank $(k-1)$.
\end{thm}

\medskip\noindent
{\bf Proof:} Let $P\subseteq \mathbb{R}^k$ be the closed, convex cone
\begin{equation}\label{eqn:conus-P}
P = \{(x_1,\dots, x_k)\in \mathbb{R}^{k} \mid (\forall i)\,  x_i\geqslant 0 \mbox{ {\rm and} } (\exists x_0\geqslant 0)\, x_0n = x_1l_1+\dots + x_kl_k \}\, .
\end{equation}
Let $W = P\cap \mathbb{Z}^{k}$ be the abelian semigroup of all lattice points in $P$. Obviously $\mathbbm{1} = (1,\dots, 1)$ is in $W$. Let $\mathbb{Z}_{\geqslant 0} \, \mathbbm{1} = \{m\mathbbm{1} \vert\, m\in \mathbb{Z}_{\geqslant 0}\} \subset W$ be  the subsemigroup of $W$ generated by the vector $\mathbbm{1}$. Then, as a consequence of Proposition \ref{prop:monomial-cycles}, there is an isomorphism of abelian groups
\[
    H_1(Tonn^{n,k}(L); \mathbb{Z}) \cong  W/(\mathbb{Z}_{\geqslant 0} \, \mathbbm{1}) \, .
\]
It is not difficult to see that the map  $p : P \rightarrow H_0^L$, which sends $x\in P$ to $p(x) = x-x_0\mathbbm{1}$ (see (\ref{eqn:conus-P})), induces an isomorphism
$W/(\mathbb{Z}_{\geqslant 0} \, \mathbbm{1}) \longrightarrow H^L_0\cap \mathbb{Z}^k$.   For example it induces an epimorphism since for each lattice point $y\in H^L_0\cap \mathbb{Z}^k$ the vector $y + m\mathbbm{1}$ is in $W$ for a sufficiently large positive integer $m$.

Let $\lambda := l_1l_2\cdots l_k$ and $\lambda_i := \lambda/l_i$. The vectors  $z^{(j)} = (z^{(j)}_1,\dots, z^{(j)}_k)$, where $z^{(j)}_j = n\lambda_j$ and $z^{(j)}_i = 0$ for $i\neq j$, clearly belong to $W$. The corresponding vectors $\{\hat{z}^{(j)} := z^{(j)} - \lambda \mathbbm{1}\}\in H_0^L \cap \mathbb{Z}^k$ span the hyperplane $H_0^L$. From here we deduce that  ${\rm rank}(H_0^L \cap \mathbb{Z}^k) = k-1$.
\hfill $\square$

\medskip
We conclude this section by some observations about the cycles   $c_i$, introduced in Definition \ref{def:cycles}.

Suppose that $l_i$ is not relatively prime to $n$, say $n = qd$ and $l_i = pd$, where $d\geqslant 2$ and $p$ and $q$ are relatively prime. In this case the cycle $c_i$ can be decomposed as a sum of $d$ (irreducible) cycles, each of length $q$. The following proposition claims that all these cycles determine the same homology class.

\begin{prop}
  All cycles of the form ${}_vE_i^q$ are homologous.
\end{prop}

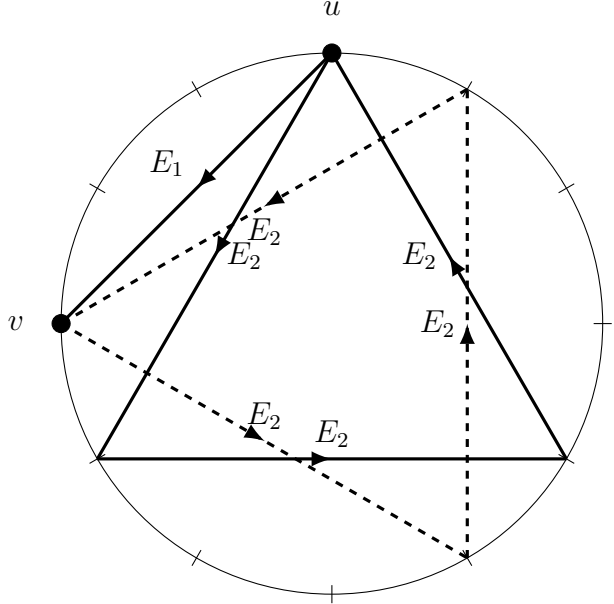
\begin{figure}[htb]
	\centering
	\begin{tikzpicture}[scale=1.2]
	\newcommand{\ang}{360/12}	% angle
	\newcommand{\R}{3}	% radius
	\draw (0,0) circle(3);
        
    \foreach \x in {0,...,8,9,10,11} {
     	\draw ({\x*\ang}:{\R-.1}) -- ({\x*\ang}:{\R+.1});
     	%\node () at ({\x*\ang}:{\R+.5}) {\x};
  	}
        
    \draw[very thick,->-] ({3*\ang}:\R) to node[right] {$E_2$} ({7*\ang}:{\R}) ;
    \draw[very thick,->-] ({7*\ang}:\R) to node[above] {$E_2$} ({11*\ang}:{\R}) ;
    \draw[very thick,->-] ({11*\ang}:\R) to node[left] {$E_2$} ({3*\ang}:{\R}) ;
        
    \draw[fill=black] ({6*\ang}:\R) circle(0.1);
    \node () at ({6*\ang}:{\R+.5}) {$v$};
    \draw[fill=black] ({3*\ang}:\R) circle(0.1);
    \node () at ({3*\ang}:{\R+.5}) {$u$};
        
    \draw[very thick,->-] ({3*\ang}:\R) to node[above left] {$E_1$} ({6*\ang}:{\R}) ;
    \draw[very thick,dashed,->-] ({6*\ang}:\R) to node[above] {$E_2$} ({10*\ang}:{\R});
    \draw[very thick,dashed,->-] ({10*\ang}:\R) to node[left] {$E_2$} ({2*\ang}:{\R});
    \draw[very thick,dashed,->-] ({2*\ang}:\R) to node[below] {$E_2$} ({6*\ang}:{\R});
\end{tikzpicture}
    \caption{${}_uE_1 E_2 ^3 E_1 ^{-1} ={}_uE_2 ^3$}
    \label{fig-3}
\end{figure}

\medskip\noindent
{\bf Proof:} Let $u$ be another base point such that $v = u+l_j$. Then the cycle  ${}_uE_jE_i^qE_j^{-1}$ is clearly homologous to the cycle ${}_vE_i^q$. On the other hand
\[
  {}_uE_jE_i^qE_j^{-1} = {}_uE_i^qE_jE_j^{-1} = {}_uE_i^q \, .
\]
By iterating this argument we see that ${}_xE_i^q$ and ${}_{x+z}E_i^q$ are homologous for any integer $z$ which can be written in the form $z = p_0n+p_1l_1+\dots + p_kl_k$. Since the vector $L$ is reduced we see that $z=1$ for some choice of parameters $p_0, p_1, \dots, p_k$, which completes the proof of the proposition. \hfill $\square$

\section{Proof of Theorem \ref{thm:main}}
\label{sec:Proof-main}

We already know that the fundamental group of the generalized tonnetz $T=Tonn^{n,k}(L)$ is free abelian of rank $k-1$ (Theorem \ref{thm:main-homology}). For the continuation of the proof of Theorem \ref{thm:main}, a natural step would be to show that  $Tonn^{n,k}(L)$ is an aspherical manifold, in the sense that its all  higher homotopy groups are trivial. Note however that asphericity alone is not sufficient to guarantee that such a manifold is covered by an euclidean space, see \cite{Davis} for examples.

\medskip
 We offer a direct proof that the universal covering $\widetilde{T}$ of a generalized Tonnetz $T=Tonn^{n,k}(L)$ is homeomorphic to $\mathbb{R}^{k-1}$ (with a lattice $\Lambda_L\subset \mathbb{R}^{k-1}$ of rank $(k-1)$ as a group of deck transformations) which implies that $T$ is homeomorphic to a $(k-1)$-dimensional torus.   To this end we construct a {\em discrete Abel-Jacobi map}
 \[
      \Omega : \widetilde{T} \longrightarrow \mathcal{D}_ \Lambda
 \]
 where $H_0 = \{x\in \mathbb{R}^{k} \mid x_1+\dots + x_k = 0\}$,   $\Lambda := H_0\cap \mathbb{Z}^k$ is a lattice isomorphic to the {\em permutohedral lattice} $\mathbb{A}_{k-1}^\ast$ and $\mathcal{D}_\Lambda$ is the associated {\em Delone triangulation}.

 \subsection{Simplicial universal covering}
 \label{sec:univ-covering}

 It is well known that each finite simplicial complex $K$ admits an universal covering $p : \widetilde{K} \rightarrow K$ in the simplicial category.

 \medskip
 By a classical construction, see Seifert-Threlfall \cite{ST-1980}, the vertices of $\widetilde{K}$ are (combinatorial)  homotopy classes of (simplicial) edge-paths $\alpha = \alpha_{\, x} = {}_{v_0}\alpha_{\, x}$ in $K$, connecting the base-point $v_0$ with a (variable) vertex $x\in K$.

 \medskip
 By definition a simplex in $\widetilde{K}$ is a collection of edge-paths   $\{{\alpha_i}_{\, x_i}\}_{i=0}^d$  (or rather their homotopy classes) such that the end-points form a simplex $\tau = \{x_0,\dots, x_d\}$ in $K$ and for each $i\neq j$
 the edge-paths ${\alpha_i}_{\, x_i}$ and ${\alpha_j}_{\, x_j}$  are {\em neighbours} in the sense of the following definition.

\begin{defin}
Two edge paths $\alpha_{\, x}$ and $\beta_{\, y}$ are {\em neighbors} if $\{x,y\}$ is an edge $e\in K$ and the edge-path $\alpha_{\, x}e\beta_{\, y}^{-1}$
is a homotopically trivial loop based at $v_0$.
 \end{defin}
We emphasize that the homotopy always refers to combinatorial homotopy. In particular two edge-paths ${}_a\alpha_b$ and ${}_a\beta_b$ are homotopic means that one can be obtained from the other by a sequence of elementary modifications (moves) which replace one side of a triangle by the remaining two sides (or vice versa).

 \subsection{Discrete Abel-Jacobi map}
\label{sec:Abel-Jacobi}

Following the notation from Section \ref{sec:homology}, each edge-path $\alpha = {}_{v_0}\alpha_x$, which emanates from the base point $v_0$ and ends at a vertex $x$, is homotopic (rel the end-points) to an edge-path of the form
\begin{equation}\label{eqn:edge-paths}
      \alpha = {}_{v_0}E_1^{p_1}E_2^{p_2}\cdots {E_k^{p_k}}_{ x} =  E_1^{p_1}E_2^{p_2}\cdots E_k^{p_k}
\end{equation}
where $p_i\geqslant 0$ for each $i\in [k]$.

\medskip
The canonical cocycles $\omega_i$, introduced in Section \ref{sec:canonical}, together define  a vector valued  $1$-cocycle $\omega = (\omega_1,\omega_2,\dots,\omega_k)$   on the generalized tonnetz $Tonn^{n,k}(L)$ which in light of (\ref{eqn:omega-zavisnost})
takes values in the subspace $H_0 = \{y\in\mathbb{R}^k \mid y_1+\dots + y_k = 0 \}\subset \mathbb{R}^k$. More precisely, the cocycle $\omega$ takes values in the lattice
\[
      \Lambda = \{x\in \mathbb{Z}^k \mid \sum_{i=1}^{k}x_i = 0 \mbox{ {\rm and all $x_i$ are in the same (mod $k$) congruence class}}\} \, .
\]
This is one of incarnations of the lattice of type $\mathbb{A}^\ast_{k-1}$ (the dual of the root lattice $\mathbb{A}_{k-1}$), which can be also described as the projection on $H_0$ of the $k$-fold dilatation $k\mathbb{Z}^k$  of the cubical lattice $\mathbb{Z}^k$, along the main diagonal $D = \{(t,t,\dots, t)\}_{t\in \mathbb{R}}\subset \mathbb{R}^k$.

\medskip
In the sequel we will need a more precise description of {\em Delone cells (simplices)} of this lattice.
Following \cite{conway2013sphere}, for $i=0, 1,\dots, k$ let $[i] = (j^i, (-i)^j)= (j,\dots, j, -i, \dots, -i)\in \mathbb{R}^k$, where $i+j=k$ and in the vector $[i]$ there are $i$ occurrences of $j$ (respectively $j$ occurrences of $-i$). Similarly if $\pi\in S_k$ is a permutation then $[i]_\pi = \pi([i])$ is obtained from $[i]$ by permuting the coordinates.
\begin{rem}\label{rem:a-notation}{\rm
Denote $a_i := \omega(E_i)$. Then $a_1+a_2+\dots+ a_k=0$ and $\{a_i\}_{i\neq j}$ is a basis of the lattice $\Lambda$ for each $j\in [k]$. For $I = \{i_1, \dots, i_r\}\subseteq [k]$ let $a_I:= \sum_{j\in I} a_i$ and $E_I = E_{i_1}\cdots E_{i_r}$. It is easily checked that the vector $[i]$ (in the traditional notation \cite{conway2013sphere}) is the same as the vector $a_{[i]} = \omega(E_1\cdots E_i)= \omega(E_{[i]})$.
}
\end{rem}

\begin{prop}\label{prop:Delone}{\rm (\cite[Theorem 4.5]{BAD-2013}, \cite[Chapters 4 and 21]{conway2013sphere})}
The Delone cells of the lattice $\Lambda \cong \mathbb{A}^\ast_{k-1}$ are $(k-1)$-simplices, which are related via permutations of coordinates and translation to the canonical simplex whose vertices are $[0]=[k], [1], \dots, [k-1]$.
\end{prop}

\begin{rem}\label{rem:a-root-lattice} {\rm
In the notation of Remark \ref{rem:a-notation} the vertices of a Delone cell are the vectors $a_{[j]}=[j]$. Each element of the lattice $\mathbb{A}_{k-1}^\ast$ has a representation $z = p_1a_1+p_2a_2+\dots + p_ka_k$ where $p_j\in \mathbb{Z}$. (This representation is unique if $p_1+\dots+ p_k = 0$.) Let $Star_{\mathcal{D}_\Lambda}(z)$ be the union of all Delone cells which have $z$ as a vertex (the star of $z$ in the Delone complex $\mathcal{D}_\Lambda$). Then the vertices of simplices in $Star_{\mathcal{D}_\Lambda}(z)$ are the vectors $z+a_I$ for all subsets $I\subseteq [k]$ and each simplex $\sigma  \in Star_{\mathcal{D}_\Lambda}(z)$ is of the form
\[
   \sigma = {\rm Conv}\{z+a_{I_1}, z+a_{I_2}, \dots, z +a_{I_s}\}
\]
where $I_1\subset I_2\subset \dots \subset I_s$.
}
\end{rem}

\medskip
The vector valued cocycle $\omega$ can be extended to edge-paths (\ref{eqn:edge-paths}) by the formula
\begin{equation}\label{eqn:omega-paths}
       \omega(\alpha)  = \omega(E_1^{p_1}E_2^{p_2} \cdots E_k^{p_k}) =  \sum_{i=1}^{k}  p_i\omega(E_i) = \sum_{i=1}^{k}  p_i a_i \in \Lambda \, .
\end{equation}
Since $\omega$ is a cocycle, $\omega({}_a\alpha_b) = \omega({}_a\beta_b)$ for each two homotopic edge-paths (with the same end-points).
The following proposition says that $\omega$ takes different values on non-homologous cycles.

\begin{prop}\label{prop:mono}
  The map $\omega$, described by the formula (\ref{eqn:omega-paths}), induces a monomorphism
  \[
    \widehat{\omega} : H_1(Tonn^{n,k}(L); \mathbb{Z}) \longrightarrow \Lambda \, .
  \]
\end{prop}
\medskip\noindent
{\bf Proof:} It is sufficient to show that if $\alpha = E_1^{p_1}E_2^{p_2}\cdots E_{k}^{p_k}$ is a loop such that $\widehat{\omega}(\alpha) = 0$ then $\alpha$ is trivial in the homology group. However, if $\omega_i(\alpha) = 0$ for each $i\in [k]$ then
\[
     k\theta_{i,j}(\alpha) = \omega_i(\alpha) - \omega_j(\alpha) =0
\]
for each $i\neq j$. In turn $\theta_{i,j}(\alpha) = 0$ for each $i\neq j$ which implies $p_1=p_2=\dots =p_k$ and, in light of Proposition \ref{prop:monomial-cycles}, $\alpha$ is a trivial cycle.
\hfill $\square$

\begin{prop}\label{prop:same-ends}
Suppose that two edge paths $\alpha = {}_a\alpha_x$ and $\beta = {}_a\beta_y$,  which share the same initial point $a$, satisfy the equality $\omega(\alpha) = \omega(\beta)$. Then $x=y$, i.e.\ they have the same end-point as well.
\end{prop}

\medskip\noindent
{\bf Proof:} Suppose that $\alpha = {{}_a}E_1^{p_1}E_2^{p_2}\cdots E_k^{p_k}{_x}$ and $\beta = {{}_a}E_1^{q_1}E_2^{q_2}\cdots E_{k}^{q_k}{_y}$. By assumption
\[
\omega(\alpha) = \sum_{i=1}^k p_i a_i = \sum_{i=1}^k q_i a_i = \omega(\beta) \, .
\]
It follows that $p_i - q_i$ does not depend on $i$ and the result follows. \hfill $\square$

\medskip
The following proposition refines Proposition \ref{prop:mono}. It says that a vector $b \in {\rm Image}(\widehat{\omega})$ cannot be ``very short'' (unless it is zero).

\begin{prop}\label{prop:refine}
Suppose that $\xi = E_1^{p_1}E_2^{p_2}\cdots E_{k}^{p_k}$ is cycle representing a non-trivial homology class in $H_1(Tonn^{n,k}(L))$, where $p_i\geqslant 0$ for each $i\in [k]$. Then $p_i\geqslant 3$ for some $i$. Moreover, its image ${\omega}(\xi)\in
 \Lambda$ in the lattice $\Lambda$ cannot be expressed as a difference $a_I - a_J$ of vectors described in Remark \ref{rem:a-notation}, where $I$ and $J$ are subsets of $[k]$.
\end{prop}

\medskip\noindent
{\bf Proof:} Since $\xi_0 = E_1E_2\cdots E_{k}$ is a trivial cycle, by factoring out from $\xi$ the power $(\xi_0)^\nu$, where $\nu:= {\rm min}\{p_j\}_{j=1}^k$, we can assume that $p_j=0$ for some $j$.

Since $\xi$ is a cycle we know that $p_1l_1+\dots+p_kl_k = p_0n$ is divisible by $n$. If $p_0 \geqslant 2$ then $p_i\geqslant 3$ for some $i\in [k]$ and we are done. Otherwise $p_0=1$ and $0\leqslant p_i\leqslant 2$ for each $i$. This is not possible since the equality
\[
     (p_1l_1+\dots+p_kl_k) - (l_1+\dots + l_k) = n-n=0
\]
would contradict the genericity of the vector $L$.  As an immediate consequence we see that the equality
${\omega}(\xi) = a_I - a_J = a_I - (a_{[k]} - a_{J^c}) = a_I+a_{J^c}$
is not possible.  \hfill $\square$

\medskip
Formula (\ref{eqn:omega-paths})  can be used for the definition of a simplicial map $\Omega : \widetilde{T} \longrightarrow  \mathcal{D}_\Lambda $ where $\mathcal{D}_\Lambda$ is the {\em Delone triangulation} of the $(k-1)$-dimensional, affine space $H_0\subset \mathbb{Z}^k$, associated to the lattice $\Lambda$.

\medskip
More explicitly, if $\widetilde{\tau} = \{{\alpha_i}_{\, x_i}\}_{i=0}^d$ is a simplex in $\widetilde{T}$, then $\Omega(\widetilde{\tau}) = \{\omega({\alpha_0}_{\, x_0}), \dots, \omega({\alpha_d}_{\, x_d})\}$.

\begin{prop}\label{prop:iso}
  The map  $\Omega : \widetilde{T} \longrightarrow  \mathcal{D}_\Lambda $ is an isomorphism of simplicial complexes.
\end{prop}

\medskip\noindent
{\bf Proof:} The map $\Omega$ is clearly an epimorphism on vertices. Indeed, if $z = p_1a_1+\dots+ p_ka_k$ is a vertex of $\mathcal{D}_\Lambda$ then $z = \Omega(E_1^{p_1}\cdots E_k^{p_k})$.

We continue by showing that $\Omega$ is a local isomorphism of simplicial complexes.  (In particular $\Omega$ is a simplicial map.)
Recall that $[i] = \omega(E_1E_2\cdots E_i) = a_{[i]}$ for each $i$ (including the case $i=0$ when we have the empty word). Similarly $[i]_\pi = \omega(E_{\pi(1)}E_{\pi(2)}\cdots E_{\pi(i)})$ for each permutation $\pi\in S_k$.

Let $\alpha = {}_{{v_0}}\alpha_x$ be an edge-path describing a vertex in $\widetilde{T}$ (connecting the base point $v_0$ with a vertex $x$ in $T$). Then in light of Proposition \ref{prop:manifold} (see also Remark \ref{rem:zonotope}) the star $Star_{\widetilde{T}}(\alpha)$ of this vertex  is the union of $k!$ simplices (one for each $\pi\in S_k$)
\[
     \widetilde{\tau}_\pi = \{\alpha, \alpha E_{\pi(1)}, \alpha E_{\pi(1)}E_{\pi(2)}, \dots , \alpha E_{\pi(1)}E_{\pi(2)} E_{\pi(k-1)}  \} \, .
\]
In light of (\ref{eqn:omega-paths}) the $\Omega$-image of this simplex is
\[
   \Omega(\widetilde{\tau}_\pi) = \{\omega(\alpha), \omega(\alpha) + [1]_{\pi}, \omega(\alpha) + [2]_{\pi}, \dots ,  \omega(\alpha) + [k-1]_{\pi} \} \, .
\]
It follows from Propositions \ref{prop:manifold} (Remark \ref{rem:zonotope}) and Proposition \ref{prop:Delone} that $\Omega$ maps bijectively the star $Star_{\widetilde{T}}(\alpha)$ of $\alpha$ in $\widetilde{T}$ to the star $Star_{\mathcal{D}_{\lambda}}(\omega(\alpha))$ of $\omega(\alpha)$ in the Delone triangulation of $H_0$.

\medskip
The map $\Omega$ is actually a covering projection. For this it is sufficient to show that for each $z \in \Lambda$ the open star $OpStar_{\mathcal{D}_\Lambda}(z) = {\rm Int}(Star_{\mathcal{D}_\Lambda}(z))$ is evenly covered by open stars  in $\widetilde{T}$. More precisely we demonstrate that the inverse image
\[
   \Omega^{-1}(OpStar_{\mathcal{D}_\Lambda}(z))= \bigcup_{\omega(\alpha) = z}  OpStar_{\widetilde{T}}(\alpha)
\]
is a disjoint union of open stars in $\widetilde{T}$.
Let  $\alpha = {}_a\alpha_x$ and $\beta = {}_a\beta_y$ be two edge-paths representing two vertices in $\widetilde{T}$.
We want to show that if $ \omega({}_{v_0}\alpha_x) = z = \omega({}_{v_0}\beta_y) $ and $Star_{\widetilde{T}}(\alpha)\cap Star_{\widetilde{T}}(\beta) \neq \emptyset$ then $\alpha$ and $\beta$ represent the same vertex in $\widetilde{T}$.

\medskip
Assume the opposite. As a consequence of Proposition \ref{prop:same-ends} we know that $x=y$, i.e.\ $\alpha$ and $\beta$ share the same end-point.  It follows that $\xi = \alpha - \beta$ is a cycle in $T$ which defines a non-trivial homology class (otherwise $\alpha$ and $\beta$ would represent the same vertex in $\widetilde{T}$).

The intersection $K = Star_{\widetilde{T}}(\alpha)\cap Star_{\widetilde{T}}(\beta)$ is a subcomplex of both stars. If this intersection is non-empty then it contains a vertex $e$ of both stars, hence $\omega(e) = \omega(\alpha) + a_I = \omega(\beta)+ a_J$ for some subsets $I$ and $J$ of $[k]$. This implies that
\[
     \Omega(\xi) = \omega(\xi) = \omega(\alpha) - \omega(\beta) = a_I - a_J = a_I + a_{J^c} \, .
\]
However, Proposition \ref{prop:refine} says that this is not possible. In other words  the cycle $\xi$ has too small image $\omega(\xi)$ for a non-trivial homology class. Hence the cycle $\alpha-\beta$ is trivial and the edge paths $\alpha$ and $\beta$  represent the same vertex in the universal cover $\widetilde{T}$.

\medskip
Finally, since $\widetilde{T}$ is connected and $\mathcal{D}_\Lambda$ is simply connected,  we conclude that the covering map $\Omega$ must be  an isomorphism of simplicial complexes.    \hfill $\square$

\bigskip\noindent
{\bf Completion of the proof of Theorem \ref{thm:main}:} The isomorphism $\Omega$ is clearly $\Gamma$-equivariant, where $\Gamma = H_1(Tonn^{n,k}(L); \mathbb{Z})$ acts on $\mathcal{D}_\Lambda$ via the monomorphism $\widehat{\omega}$ from Proposition \ref{prop:mono} (see also the formula   (\ref{eqn:omega-paths})). It immediately follows that $Tonn^{n,k}(L)$ isomorphic to the simplicial complex $\mathcal{D}_\Lambda/\Lambda_L$ where $\Lambda_L := \widehat{\omega}(\Gamma)\subset \Lambda$  is a free abelian group of rank $k-1$.  \hfill $\square$

\section{Examples and concluding remarks}

The isomorphism  $\Omega$, described in Proposition \ref{prop:iso}, can be used  for comparison of combinatorial types of different  complexes of Tonnetz type.

\medskip
Note that each automorphism of the Delone simplicial complex $\mathcal{D}_\Lambda$ induces an isometry on the ambient euclidean space $H_0\subset \mathbb{R}^k$.
If two simplicial complexes $T_1$ and $T_2$ of Tonnetz type are combinatorially isomorphic then their universal covers $\widetilde{T}_1$ and $\widetilde{T}_2$ are also combinatorially isomorphic.

\medskip
From these two observation we conclude that each Tonnetz inherits a canonical metric from the euclidean space $H_0$ which is an invariant of its combinatorial type.

\begin{figure}[htb]
    \centering
    %\hfill
    \subfigure[$Tonn^{12,3} \left(3,4,5\right)$]{\begin{tikzpicture}[scale=1.2]
    \node[draw=black,circle] (r1n1) at (0,0) {1};
    \node[draw=black,circle] (r1n2) at (1,0) {4};
    \node[draw=black,circle] (r1n3) at (2,0) {7};
    \node[draw=black,circle] (r1n4) at (3,0) {10};
    \node[draw=black,circle] (r1n5) at (4,0) {1};
    \node[draw=black,circle] (r2n1) at (1/2,{-sqrt(3)/2}) {9};
    \node[draw=black,circle] (r2n2) at (3/2,{-sqrt(3)/2}) {0};
    \node[draw=black,circle] (r2n3) at (5/2,{-sqrt(3)/2}) {3};
    \node[draw=black,circle] (r2n4) at (7/2,{-sqrt(3)/2}) {6};
    \node[draw=black,circle] (r2n5) at (9/2,{-sqrt(3)/2}) {9};
    \node[draw=black,circle] (r3n1) at (1,{-2*sqrt(3)/2}) {5};
    \node[draw=black,circle] (r3n2) at (2,{-2*sqrt(3)/2}) {8};
    \node[draw=black,circle] (r3n3) at (3,{-2*sqrt(3)/2}) {11};
    \node[draw=black,circle] (r3n4) at (4,{-2*sqrt(3)/2}) {2};
    \node[draw=black,circle] (r3n5) at (5,{-2*sqrt(3)/2}) {5};
    \node[draw=black,circle] (r4n1) at (3/2,{-3*sqrt(3)/2}) {1};
    \node[draw=black,circle] (r4n2) at (5/2,{-3*sqrt(3)/2}) {4};
    \node[draw=black,circle] (r4n3) at (7/2,{-3*sqrt(3)/2}) {7};
    \node[draw=black,circle] (r4n4) at (9/2,{-3*sqrt(3)/2}) {10};
    \node[draw=black,circle] (r4n5) at (11/2,{-3*sqrt(3)/2}) {1};
    
    \draw (r1n1) -- (r1n2) -- (r1n3) --(r1n4) -- (r1n5);
    \draw (r2n1) -- (r2n2) -- (r2n3) --(r2n4) -- (r2n5);
    \draw (r3n1) -- (r3n2) -- (r3n3) --(r3n4) -- (r3n5);
    \draw (r4n1) -- (r4n2) -- (r4n3) --(r4n4) -- (r4n5);
    
    \draw (r1n1) -- (r2n1) -- (r3n1) -- (r4n1);
    \draw (r1n2) -- (r2n2) -- (r3n2) -- (r4n2);
    \draw (r1n3) -- (r2n3) -- (r3n3) -- (r4n3);
    \draw (r1n4) -- (r2n4) -- (r3n4) -- (r4n4);
    \draw (r1n5) -- (r2n5) -- (r3n5) -- (r4n5);

    \draw (r2n1) -- (r1n2);
    \draw (r3n1) -- (r2n2) -- (r1n3);
    \draw (r4n1) -- (r3n2) -- (r2n3) -- (r1n4);
    \draw (r4n2) -- (r3n3) -- (r2n4) -- (r1n5);
    \draw (r4n3) -- (r3n4) -- (r2n5);
    \draw (r4n4) -- (r3n5);
\end{tikzpicture}}
    %\hfill
    \subfigure[$Tonn^{12,3}\left(2,3,7\right)$]{\begin{tikzpicture}[scale=1.2]
    \node[draw=black,circle] (u2n1) at (2,{2*sqrt(3)/2}) {10};
    \node[draw=black,circle] (u1n1) at (3/2,{sqrt(3)/2}) {5};
    \node[draw=black,circle] (u1n2) at (5/2,{sqrt(3)/2}) {7};
    \node[draw=black,circle] (r1n1) at (0,0) {10};
    \node[draw=black,circle] (r1n2) at (1,0) {0};
    \node[draw=black,circle] (r1n3) at (2,0) {2};
    \node[draw=black,circle] (r1n4) at (3,0) {4};
    \node[draw=black,circle] (r2n1) at (1/2,{-sqrt(3)/2}) {7};
    \node[draw=black,circle] (r2n2) at (3/2,{-sqrt(3)/2}) {9};
    \node[draw=black,circle] (r2n3) at (5/2,{-sqrt(3)/2}) {11};
    \node[draw=black,circle] (r2n4) at (7/2,{-sqrt(3)/2}) {1};
    \node[draw=black,circle] (r3n1) at (1,{-2*sqrt(3)/2}) {4};
    \node[draw=black,circle] (r3n2) at (2,{-2*sqrt(3)/2}) {6};
    \node[draw=black,circle] (r3n3) at (3,{-2*sqrt(3)/2}) {8};
    \node[draw=black,circle] (r3n4) at (4,{-2*sqrt(3)/2}) {10};
    \node[draw=black,circle] (r4n1) at (3/2,{-3*sqrt(3)/2}) {1};
    \node[draw=black,circle] (r4n2) at (5/2,{-3*sqrt(3)/2}) {3};
    \node[draw=black,circle] (r4n3) at (7/2,{-3*sqrt(3)/2}) {5};
    \node[draw=black,circle] (r5n1) at (2,{-4*sqrt(3)/2}) {10};
    \node[draw=black,circle] (r5n2) at (3,{-4*sqrt(3)/2}) {0};
    
    \draw (u1n1) -- (u1n2) -- (u2n1) -- (u1n1);
    \draw (r1n1) -- (r1n2) -- (r1n3) -- (r1n4);
    \draw (r1n2) -- (u1n1) -- (r1n3) -- (u1n2) -- (r1n4);
    \draw (r2n1) -- (r2n2) -- (r2n3) -- (r2n4);
    \draw (r3n1) -- (r3n2) -- (r3n3) -- (r3n4);
    \draw (r4n1) -- (r4n2) -- (r4n3);
    \draw (r1n1) -- (r2n1) -- (r1n2) -- (r2n2) -- (r1n3) -- (r2n3) -- (r1n4) -- (r2n4);
    \draw (r2n1) -- (r3n1) -- (r2n2) -- (r3n2) -- (r2n3) -- (r3n3) -- (r2n4) -- (r3n4);
    \draw (r3n1) -- (r4n1) -- (r3n2) -- (r4n2) -- (r3n3) -- (r4n3) -- (r3n4);
    \draw (r4n1) -- (r5n1) -- (r4n2) -- (r5n2) -- (r4n3);
    \draw (r5n1) -- (r5n2);
\end{tikzpicture}}
    %\hfill
    \caption{Combinatorially non-isomorphic complexes of Tonnetz type.}
    \label{fig:prva}
\end{figure}
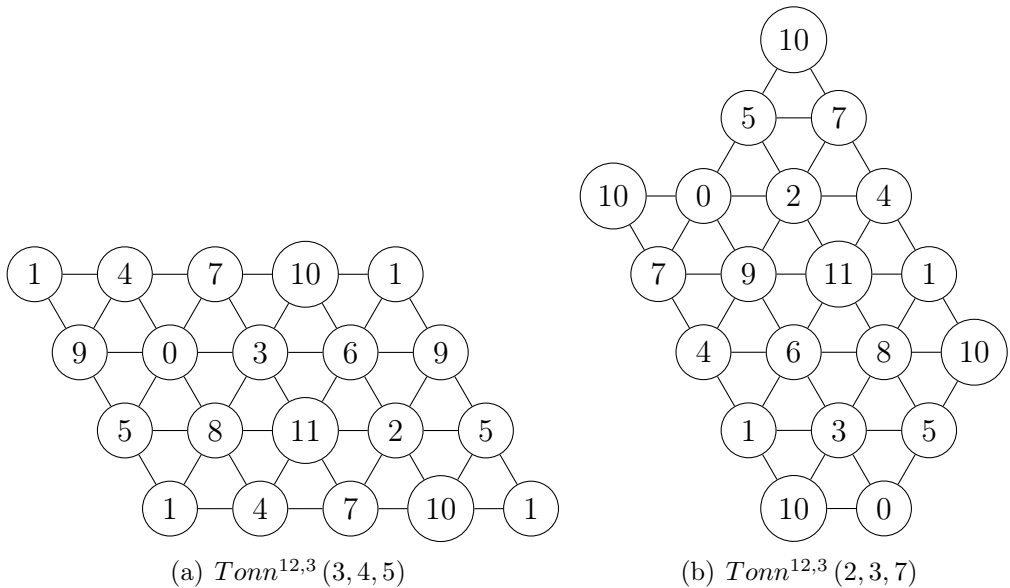

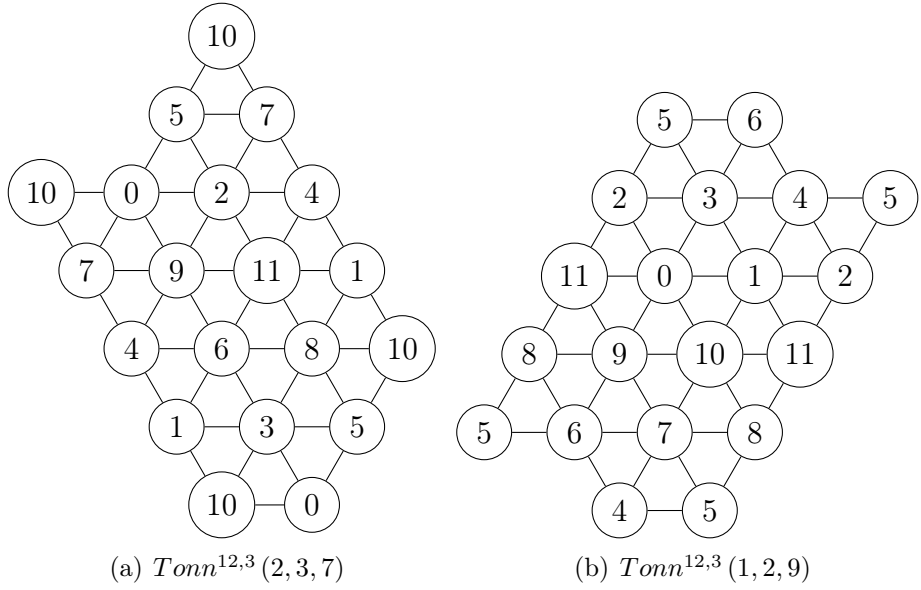
\begin{figure}[htb]
    \centering
    %\hfill
    \subfigure[$Tonn^{12,3} \left(2,3,7\right)$]{\begin{tikzpicture}[scale=1.2]
    \node[draw=black,circle] (u2n1) at (2,{2*sqrt(3)/2}) {10};
    \node[draw=black,circle] (u1n1) at (3/2,{sqrt(3)/2}) {5};
    \node[draw=black,circle] (u1n2) at (5/2,{sqrt(3)/2}) {7};
    \node[draw=black,circle] (r1n1) at (0,0) {10};
    \node[draw=black,circle] (r1n2) at (1,0) {0};
    \node[draw=black,circle] (r1n3) at (2,0) {2};
    \node[draw=black,circle] (r1n4) at (3,0) {4};
    \node[draw=black,circle] (r2n1) at (1/2,{-sqrt(3)/2}) {7};
    \node[draw=black,circle] (r2n2) at (3/2,{-sqrt(3)/2}) {9};
    \node[draw=black,circle] (r2n3) at (5/2,{-sqrt(3)/2}) {11};
    \node[draw=black,circle] (r2n4) at (7/2,{-sqrt(3)/2}) {1};
    \node[draw=black,circle] (r3n1) at (1,{-2*sqrt(3)/2}) {4};
    \node[draw=black,circle] (r3n2) at (2,{-2*sqrt(3)/2}) {6};
    \node[draw=black,circle] (r3n3) at (3,{-2*sqrt(3)/2}) {8};
    \node[draw=black,circle] (r3n4) at (4,{-2*sqrt(3)/2}) {10};
    \node[draw=black,circle] (r4n1) at (3/2,{-3*sqrt(3)/2}) {1};
    \node[draw=black,circle] (r4n2) at (5/2,{-3*sqrt(3)/2}) {3};
    \node[draw=black,circle] (r4n3) at (7/2,{-3*sqrt(3)/2}) {5};
    \node[draw=black,circle] (r5n1) at (2,{-4*sqrt(3)/2}) {10};
    \node[draw=black,circle] (r5n2) at (3,{-4*sqrt(3)/2}) {0};
    
    \draw (u1n1) -- (u1n2) -- (u2n1) -- (u1n1);
    \draw (r1n1) -- (r1n2) -- (r1n3) -- (r1n4);
    \draw (r1n2) -- (u1n1) -- (r1n3) -- (u1n2) -- (r1n4);
    \draw (r2n1) -- (r2n2) -- (r2n3) -- (r2n4);
    \draw (r3n1) -- (r3n2) -- (r3n3) -- (r3n4);
    \draw (r4n1) -- (r4n2) -- (r4n3);
    \draw (r1n1) -- (r2n1) -- (r1n2) -- (r2n2) -- (r1n3) -- (r2n3) -- (r1n4) -- (r2n4);
    \draw (r2n1) -- (r3n1) -- (r2n2) -- (r3n2) -- (r2n3) -- (r3n3) -- (r2n4) -- (r3n4);
    \draw (r3n1) -- (r4n1) -- (r3n2) -- (r4n2) -- (r3n3) -- (r4n3) -- (r3n4);
    \draw (r4n1) -- (r5n1) -- (r4n2) -- (r5n2) -- (r4n3);
    \draw (r5n1) -- (r5n2);
\end{tikzpicture}}
    %\hfill
    \subfigure[$Tonn^{12,3}\left(1,2,9\right)$]{\begin{tikzpicture}[scale=1.2]
    \node[draw=black,circle] (r1n1) at (0,0) {5};
    \node[draw=black,circle] (r1n2) at (1,0) {6};
    \node[draw=black,circle] (r2n1) at (-1/2,{-sqrt(3)/2}) {2};
    \node[draw=black,circle] (r2n2) at (1/2,{-sqrt(3)/2}) {3};
    \node[draw=black,circle] (r2n3) at (3/2,{-sqrt(3)/2}) {4};
    \node[draw=black,circle] (r2n4) at (5/2,{-sqrt(3)/2}) {5};
    \node[draw=black,circle] (r3n1) at (-1,{-2*sqrt(3)/2}) {11};
    \node[draw=black,circle] (r3n2) at (0,{-2*sqrt(3)/2}) {0};
    \node[draw=black,circle] (r3n3) at (1,{-2*sqrt(3)/2}) {1};
    \node[draw=black,circle] (r3n4) at (2,{-2*sqrt(3)/2}) {2};
    \node[draw=black,circle] (r4n1) at (-3/2,{-3*sqrt(3)/2}) {8};
    \node[draw=black,circle] (r4n2) at (-1/2,{-3*sqrt(3)/2}) {9};
    \node[draw=black,circle] (r4n3) at (1/2,{-3*sqrt(3)/2}) {10};
    \node[draw=black,circle] (r4n4) at (3/2,{-3*sqrt(3)/2}) {11};
    \node[draw=black,circle] (r5n1) at (-2,{-4*sqrt(3)/2}) {5};
    \node[draw=black,circle] (r5n2) at (-1,{-4*sqrt(3)/2}) {6};
    \node[draw=black,circle] (r5n3) at (0,{-4*sqrt(3)/2}) {7};
    \node[draw=black,circle] (r5n4) at (1,{-4*sqrt(3)/2}) {8};
    \node[draw=black,circle] (r6n3) at (-1/2,{-5*sqrt(3)/2}) {4};
    \node[draw=black,circle] (r6n4) at (1/2,{-5*sqrt(3)/2}) {5};
    
    \draw (r5n1) -- (r4n1) -- (r3n1) -- (r2n1) -- (r1n1);
    \draw (r5n2) -- (r4n2) -- (r3n2) -- (r2n2) -- (r1n2);
    \draw (r6n3) -- (r5n3) -- (r4n3) -- (r3n3) -- (r2n3);
    \draw (r6n4) -- (r5n4) -- (r4n4) -- (r3n4) -- (r2n4);
    
    \draw (r1n1) -- (r1n2);
    \draw (r2n1) -- (r2n2) -- (r2n3) -- (r2n4);
    \draw (r3n1) -- (r3n2) -- (r3n3) -- (r3n4);
    \draw (r4n1) -- (r4n2) -- (r4n3) -- (r4n4);
    \draw (r5n1) -- (r5n2) -- (r5n3) -- (r5n4);
    \draw (r6n3) -- (r6n4);
    
    \draw (r4n1) -- (r5n2) -- (r6n3);
    \draw (r3n1) -- (r4n2) -- (r5n3) -- (r6n4);
    \draw (r2n1) -- (r3n2) -- (r4n3) -- (r5n4);
    \draw (r1n1) -- (r2n2) -- (r3n3) -- (r4n4);
    \draw (r1n2) -- (r2n3) -- (r3n4);
    
\end{tikzpicture}}
    %\hfill
    \caption{Combinatorially isomorphic complexes of Tonnetz type.}
    \label{ex-3}
\end{figure}

\medskip
For illustration the classical Tonnetz, exhibited in Figure \ref{fig:prva}, is non-isometric (and therefore combinatorially non-isomorphic) to the ``Tonnetz'' shown in the same figure on the right.

\medskip
This can be proved by comparing the shortest closed, non-contractible geodesics (systoles) of both complexes. For example the systole on the left has the length $3$, while on the right the length is $\sqrt{7}$.

\medskip
The complex $Tonn^{12,3}(2,3,7)$ and the complex $Tonn^{12,3}(1,2,9)$ (exhibited in Figure \ref{ex-3}) are isometric, in particular have systoles of the same length. Moreover, they are combinatorially isomorphic. Indeed, if we cut out  the parallelogram $5$-$3$-$4$-$6$ from Figure \ref{ex-3}~(b) and glue it on the opposite side, we obtain a fundamental domain of the ``Tonnetz'' $Tonn^{12,3}(1,2,9)$ which, by an automorphism of the planar Delone complex  $\mathcal{D}_\Lambda$, can be mapped to the Figure \ref{ex-3}~(a).

\medskip
Figures \ref{fig:prva} and \ref{ex-3} were originally generated by lifting the triangulations from a Tonnetz $T$ to its universal cover $\widetilde{T}$. Informally speaking, they are obtained by gradually unfolding the complex $T$ in the plane until the picture becomes periodic.

Results from Section \ref{sec:homology} and \ref{sec:Proof-main}, as summarized in the following proposition, allow us to generate these and related pictures (for an arbitrary $Tonn^{n,k}(L)$)   directly from the input length vector $L = (l_1,l_2,\dots,l_k)$.

\begin{prop}\label{prop:generate}
The lattice $\Lambda_L := \widehat{\omega}(\Gamma)\subset \Lambda$, which appears in the isomorphism \newline $Tonn^{n,k}(L) \cong \mathcal{D}_\Lambda/\Lambda_L$, has the following explicit description
\[
\Lambda_L = \{p_1a_1+\dots +p_ka_k \mid (\forall i)\,  p_i\in \mathbb{Z} \mbox{ {\rm and} } (\exists p_0\in \mathbb{Z}) \, p_1l_1+\dots +p_kl_k = p_0n\}
\]
where $a_i = \omega(E_i)\in \Lambda$ are the vectors introduced in Remark \ref{rem:a-notation}.
\end{prop}

\begin{exam}
Let us explicitly describe the group $\Lambda_L$ for  $Tonn^{12,3}(2,3,7)$.
The lattice
\[
\Gamma_L = \{(x,y,z)\in \mathbb{Z}^3 \mid 2x + 3y + 7z = 0 \}
\]
has a parametric presentation
\[
\Gamma_L = \{(x,y,z)\in \mathbb{Z}^3 \mid (\exists r,s\in \mathbb{Z}) \mid\, x = -3r-5s, y = 2r+s, z = s \} \, .
\]
By choosing $(r,s) = (1,0)$ and $(r,s) = (0,1)$ we obtain that $\{(-3, 2, 0), (-5,1,1)\}$ is a basis for $\Gamma_L$.
It follows that the corresponding generators of the lattice $\Lambda_L$ are $b_1 = -3a_1+2a_2$ and $b_2 = -5a_1+a_2+a_3 = -6a_1$. By interpreting (in Figure \ref{fig:prva}) $a_1$ and $a_2$ as the vectors connecting the vertex (labeled by) $10$ by the neighbouring vertices $0$ and $1$, we easily check that vectors $b_1$ and $b_2$ preserve the labeling of this lattice. They actually generate this lattice since $\langle b_1, b_2\rangle $ is  a sublattice of $\Lambda$ of index $12$.
\end{exam}

\subsection{Irrational Tonnetz}

It is natural to extend the definition of the tonnetz to the case $n = +\infty$, interpreted as the limit case when $k$ is fixed and  $n$ approaches infinity. Informally, vertices are points on a circle $C$ with circumference $1$  while simplices are finite subsets $I\subset C$ which are $L$-admissible in the sense of the following definition.

\begin{defin}\label{def:Tonn-infinity}
Let $C = \mathbb{R}/\mathbb{Z} = [0,1]/\langle 0 \simeq 1\rangle$ be a circle with induced group structure and the corresponding  (circular) order.
Suppose that $L = (l_1,l_2,\dots, l_k)$ is a collection of positive real numbers such that:
\begin{enumerate}
    \item[{\rm (1)}] The numbers $l_i$ add up to one, $l_1+l_2+\dots+ l_k = 1$.
    \item[{\rm (2)}] $L$ is irrational in the sense that $p_1l_1+\dots+p_kl_k\neq 0$ for each $p\in \mathbb{Z}^k\setminus \{0\}$.
\end{enumerate}
A subset $I\subset C$ is $L$-admissible if there exists $x\in C$ and a permutation $\pi\in S_k$ such that
\begin{align}\label{align:tonnetz-max-simplex-infty}
    \Delta(x;\sigma)=\{x,x+l_{\sigma(1)},x+l_{\sigma(1)}+l_{\sigma(2)},\ldots,x+l_{\sigma(1)}+\ldots+l_{\sigma(k-1)}\} \, .
\end{align}
Full (irrational) tonnetz $F$-${Tonn}^{\infty, k}(L)$   is the $(k-1)$-dimensional simplicial complex of all $L$-admissible subsets of $C$. The irrational tonnetz $Tonn^{\infty, k}(L)$ is a connected component of the  full Tonnetz  $F$-${Tonn}^{\infty, k}(L)$.
\end{defin}

\begin{rem}{\rm  A rotation of the circle $C$ induces an automorphism of the full tonnetz $F$-${Tonn}^{\infty, k}(L)$. Moreover, the group $C$ acts transitively on its connected components. It follows that all connected components of the irrational tonnetz are isomorphic.

Note that the condition (2) in Definition \ref{def:Tonn-infinity} guarantees
that the length vector $L$ is generic in the sense that numbers $l_i$ satisfy the condition (\ref{eqn:generic}).
}
\end{rem}

\begin{thm}\label{thm:main-irrational}
Irrational tonnetz $Tonn^{\infty, k}(L)$ is isomorphic to the Delone triangulation of the vector space $\mathbb{R}^{k-1}$ associated to the permutohedral lattice $\mathbb{A}^\ast_{k-1}$,
\[
     Tonn^{\infty,k}(L) \cong \mathcal{D}_\Lambda \, .
\]
\end{thm}

\medskip\noindent
{\bf Proof:} Many concepts introduced in Section \ref{sec:canonical}, such as atomic $1$-simplices $E_i$, cocycles $\theta_{i,j}$, canonical cocycles $\omega_i$ etc., preserve they meaning in the case of the infinite tonnetz $Tonn^{\infty,k}(L)$. An exception are canonical cycles $c_i$ whose existence is ruled out by the condition (2) from Definition \ref{def:Tonn-infinity}. Proposition \ref{prop:abelian} still holds with essentially the same proof so the fundamental group of the infinite tonnetz is always abelian. Let us show that it is actually a trivial group.

As before each $1$-chain has a representation $X = E_1^{p_1}E_2^{p_2}\cdots E_k^{p_k}$. If this is a cycle  (with winding number $p_0$) then $p_1l_1+\dots+ p_kl_k = p_0$. In light of the condition (2) (Definition \ref{def:Tonn-infinity}) this is possible only if $p_1=\cdots = p_k = p_0$ in which case $X$ is a boundary.

The end of the proof follows closely the idea of the proof of Proposition \ref{prop:iso}. The isomorphism  $\Omega :  Tonn^{\infty,k}(L) \rightarrow \mathcal{D}_\Lambda$ is again defined with the aid of formula (\ref{eqn:omega-paths}).  \hfill $\square$

\subsection{Genericity condition (\ref{eqn:generic})}

The genericity condition (\ref{eqn:generic}) plays a central role in many  arguments. For illustration a generalized tonnetz may not be a manifold without this condition, as visible from the classification of all $2$-dimensional  (not necessarily generic) complexes of Tonnetz type, see \cite[Section 6]{MC-2011}.

\medskip
The smallest examples $(k\geqslant 3)$ of length vectors which are generic are:

\medskip
$(1,2,4)$ for $(n,k) = (7,3)$

$(1,2,5)$  for $(n,k) = (8,3)$

$(1,2,6), (2,3,4)$  for $(n,k) = (9,3)$

$(1,2,7), (1,3,6)$  for $(n,k) = (10,3)$

$(1,2,8), (1,3,7), (1,4,6), (2,3,6), (2,4,5)$  for $(n,k) = (11,3)$

$(1,2,9), (1,3,8), (1,4,7), (2,3,7), (3,4,5)$  for $(n,k) = (12,3)$, etc.

\medskip
It is not difficult to construct examples of families of generic vectors as illustrated by $(1,q,q^2,\dots, q^{k-1})$ for $q\geqslant 2$.

\medskip
When $k$ is fixed, asymptotically (when $n \rightarrow\infty$) almost all vectors are generic.
This can be deduced by observing that generic vectors are positive integer vectors in a simplex with vertices $ne_i \, (i=1,\dots, k)$ outside the union of the hyperplane arrangement  $\mathcal{H}_k = \{H_{I,J}\}$ where for two disjoint, non-empty  subspaces $I,J\subset [k]$
\[
             H_{I,J} = \{x\in \mathbb{R}^k \mid \, \sum_{i\in I} x_i = \sum_{j\in J} x_j \}  \, .
\]

\subsection{Other generalizations of the Tonnetz}

There are other generalizations of  the classical Tonnetz, see for example  \cite{MC-2011, CFS, DT-2012} or \cite{tymoczko_geometry_2006}. The authors of these papers usually  put more emphasis on combinatorial and geometric aspects of the musical theory and see mathematics primarily as a useful tool. These papers do not overlap with our exposition with an exception of  \cite{MC-2011} where the author introduced and studied the complexes of Tonnetz type  in the case $k=3$ (without the genericity condition (\ref{eqn:generic})). In particular our Proposition \ref{prop:k=3}    is included in \cite[Theorem 23]{MC-2011}. Moreover the author provides the list of all $2$-dimensional complexes which in the non-generic case  can arise as Tonnetz-type complexes.

\bigskip\noindent
{\bf Acknowledgements:} We would like to acknowledge valuable remarks and kind suggestions of the anonymous referee which helped us improve the presentation of results in the paper.

%\printbibliography

\end{document}